\documentclass[amstex,12pt, amssymb]{article}

\usepackage{mathtext}
\usepackage[cp1251]{inputenc}
\usepackage[T2A]{fontenc}
\usepackage[dvips]{graphicx}
\usepackage{amsmath}
\usepackage{amssymb}
\usepackage{amsxtra}
\usepackage{latexsym}
\usepackage{ifthen}

\newcounter{lemma}[section]

\newcounter{corol}[section]

\newcounter{rem}[section]

\newcounter{theo}[section]

\newcounter{propo}[section]

\numberwithin{equation}{section}

\begin{document}

\markboth{\centerline{E. SEVOST'YANOV}} {\centerline{ABOUT ONE
MODULUS...}}

\def\cc{\setcounter{equation}{0}
\setcounter{figure}{0}\setcounter{table}{0}}

\overfullrule=0pt

%\normalsize\large

\author{{E. SEVOST'YANOV}\\}

\title{
{\bf ABOUT ONE MODULUS INEQUALITY OF THE V\"{A}IS\"{A}L\"{A} TYPE }}

\date{\today}
\maketitle

%\large
\begin{abstract} The present paper is devoted to the study of
space mappings, which are more general than quasiregular. The
analogue of the known V\"{a}is\"{a}l\"{a} inequality for the special
classes of curves and mappings was proved for the open, discrete,
differentiable a.e. mappings having $N,$ $N^{-1}$ and $ACP^{-1}$
properties.

\end{abstract}


\section{Introduction}

Here are some definitions. Everywhere below, $D$ is a domain in
${\Bbb R}^n,$ $n\ge 2,$ $m$ is the Lebesgue measure in ${\Bbb R}^n,$
$m_1$ is the linear Lebesgue measure on ${\Bbb R}.$   A mapping
$f:D\rightarrow {\Bbb R}^n$ is said to be a {\it discrete} if the
pre-image $f^{-1}(y)$ of any point $y\,\in\,{\Bbb R}^n$ consists of
isolated points, and an {\it open} if the image of any open set
$U\subset D$ is open in ${\Bbb R}^n.$ The notation $f:D\rightarrow
{\Bbb R}^n$ assumes that $f$ is continuous.

We write $f\in W^{1,n}_{loc}(D),$ iff all of the coordinate
functions $f_j,$ $f=(f_1,\ldots,f_n),$ have the partitional
derivatives which are locally integrable in the degree $n$ in $D.$
Recall that a mapping $f:D\rightarrow {\Bbb R}^n$ is said to be {\it
a mapping with bounded distortion}, if the following conditions
hold:

\noindent 1) $f\in W_{loc}^{1,n},$ 2) a Jacobian $J(x,f):={\rm
det\,}f^{\,\prime}(x) $ of the mapping $f$ at the point $x\in D$
preserves the sign almost everywhere in $D,$ 3) $\Vert
f^{\,\prime}(x) \Vert^n \le K \cdot |J(x,f)|$ at a.e. $x\in D$ and
some constant $K<\infty,$ where $\Vert
f^{\,\prime}(x)\Vert:=\sup\limits_{h\in {\Bbb R}^n:
|h|=1}|f^{\,\prime}(x)h|,$
see $\S\, 3.$ I in \cite{Re}, or definition 2.1 of the  2. I in
\cite{Ri}.

A curve $\gamma$ in ${\Bbb R}^n$ is a continuous mapping $\gamma
:\Delta\rightarrow{\Bbb R}^n$ where $\Delta$ is an interval in
${\Bbb R} .$ Its locus $\gamma(\Delta)$ is denoted by $|\gamma|.$
Given a family of curves $\Gamma$ in ${\Bbb R}^n ,$ a Borel function
$\rho:{\Bbb R}^n \rightarrow [0,\infty]$ is called {\it admissible}
for $\Gamma ,$ abbr. $\rho \in {\rm adm}\, \Gamma ,$ if curvilinear
integral of the first type $\int\limits_{\gamma} \rho(x)|dx|$
satisfies the condition $$\int\limits_{\gamma} \rho(x)|dx| \ge 1$$
for each (locally rectifiable) $\gamma\in\Gamma.$ Given $p\ge 1,$
the {\it $p$--modulus} $M_p(\Gamma )$ of $\Gamma$ is defined as
$$M_p(\Gamma) =\inf\limits_{ \rho \in {\rm adm}\, \Gamma}
\int\limits_{{\Bbb R}^n} \rho^p(x) dm(x)$$ interpreted as $+\infty$
if ${\rm adm}\, \Gamma = \varnothing .$ The properties of it are
analogous to the properties of the measure of Lebesgue $m$ in ${\Bbb
R}^n.$ Namely, $M_p(\varnothing)=0,$ $M_p(\Gamma_1)\le
M_p(\Gamma_2)$ whenever $\Gamma_1\subset\Gamma_2$ and
$M_p\left(\bigcup\limits_{i=1}^{\infty}\Gamma_i\right)\le
\sum\limits_{i=1}^{\infty}M_p(\Gamma_i),$ see Theorem 6.2 in
\cite{Va$_2$}.

Set $l\left(f^{\,\prime}(x)\right):=\inf\limits_{h\in {\Bbb R}^n:
|h|=1}|f^{\,\prime}(x)h|.$ Recall that {\it inner dilatation of the
order $p$} of the mapping $f$ at a point $x$ is defined as
$$K_{I, p}\,(x,f)\quad =\quad= \left\{
\begin{array}{rr}
\frac{|J(x,f)|}{{l\left(f^{\,\prime}(x)\right)}^p}, & J(x,f)\ne 0,\\
1,  &  f^{\,\prime}(x)=0, \\
\infty, & {\rm otherwise}
\end{array}
\right.\,.$$
Let $K_{I, p}(f)={\rm ess}\sup\limits_{x\in D}K_{I,p}(x, f).$ Note
that $K_{I, n}(f)<K^{n-1}<\infty$ whenever $f$ is a mapping with
bounded distortion,  see (2.7) and (2.8) of 2. 1. I in \cite{Re}.

\medskip
Suppose that $\alpha$ and $\beta$ are curves in ${\Bbb R}^n.$ Then a
notation $\alpha\subset\beta$ denotes that $\alpha$ is a subpath of
$\beta.$ In what follows, $I$ denotes an open, a closed or a
semi--open interval on the real axes. The following definition can
be found in the section 5 of Ch. II in \cite{Ri}.

\medskip
Let $I=[a,b].$ Given a rectifiable path $\gamma:I\rightarrow {\Bbb
R}^n$ we define a length function $l_{\gamma}(t)$ by the rule
$l_{\gamma}(t)=S\left(\gamma, [a,t]\right),$ where $S(\gamma,
[a,t])$ is a length of the path $\gamma|_{[a, t]}.$ Let
$\alpha:[a,b]\rightarrow {\Bbb R}^n$ be a rectifiable curve in
${\Bbb R}^n,$ $n\ge 2,$ and $l(\alpha)$ be its length. A {\it normal
representation} $\alpha^0$ of $\alpha$ is defined as a curve
$\alpha^0:[0, l(\alpha)]\rightarrow {\Bbb R}^n$ which is can be got
from $\alpha$ by change of parameter such that
$\alpha(t)=\alpha^0\left(S\left(\alpha, [a, t]\right)\right).$

\medskip
Let $f:D\rightarrow {\Bbb R}^n$ be such that $f^{\,-1}(y)$ does not
contain a nondegenerate curve for any $y\in {\Bbb R}^n,$
$\beta:I_0\rightarrow {\Bbb R}^n$ be a closed rectifiable curve and
$\alpha:I\rightarrow D$ such that $f\circ \alpha\subset \beta.$ If
the length function $l_{\beta}:I_0\rightarrow [0, l(\beta)]$ is a
constant on $J\subset I,$ then $\beta$ is a constant on $J$ and
consequently a curve $\alpha$ to be a constant on $J.$ Thus, there
exists a unique function $\alpha^{\,*}:l_\beta(I)\rightarrow D$ such
that $\alpha=\alpha^{\,*}\circ (l_\beta|_I).$ We say that
$\alpha^{\,*}$ to be a {\it $f$--re\-pre\-se\-n\-ta\-tion of
$\alpha$ by the respect to $\beta$ } if $\beta=f\circ\alpha.$

\medskip
If $\alpha:[a, b]\rightarrow D$ is a closed curve, we say that $f$
winds $\alpha$ $m$ times around itself if $f\circ \alpha=\beta$ is
rectifiable and if the following condition is satisfied: Let
$\beta^0:[0, c]\rightarrow {\Bbb R}^n$ be the normal representation
of $\beta,$ let $\alpha^{\,*}:[0, c]\rightarrow D$ be
$f$--rep\-re\-sen\-ta\-tion of $\alpha$ with respect to $\beta,$ and
let $h=c/m.$ Then $\beta^0(t+jh)=\beta^0(t)$ and
$\alpha^{\,*}(t+jh)\ne \alpha^{\,*}(t)$ whenever $0\le t<t+jh<c$ and
$j\in\{1,\ldots,m-1\}.$

In 1972 in the work of J.~V\"{a}is\"{a}l\"{a} was proved the
following, see e.g. Theorem 3.9 in \cite{Va$_3$}.

\medskip
\begin{theo}\label{th1A} {\sl Let $f:D\rightarrow {\Bbb R}^n$ be a non--constant
mapping with bounded distortion. Suppose that $\Gamma$ is a curve
family in $D,$ that $m$ is a positive integer, and that $f$ winds
every path of $\Gamma$ $m$ times around itself. Then
\begin{equation}\label{eq1A}
M_n(f(\Gamma))\quad\le\quad \frac{K_{I,
n}(f)}{m}\,\cdot\,M_n(\Gamma)\,.
\end{equation}}
\end{theo}
The goal of the present paper is to prove the analogue of the
\ref{th1A} for more general classes of mappings. Recall that a
mapping $f:D\rightarrow {\Bbb R}^n$ is said to have the {\it
$N$--pro\-per\-ty (of Luzin)} if $m\left(f\left(S\right)\right)=0$
whenever $m(S)=0$ for all such sets $S\subset{\Bbb R}^n.$ Similarly,
$f$ has the {\it $N^{-1}$--pro\-per\-ty} if
$m\left(f^{\,-1}(S)\right)=0$ whenever $m(S)=0.$

We say that a property $P$ holds for {\it $p$--al\-most every
(a.e.)} curves $\gamma$ in a family $\Gamma$ if the subfamily of all
curves in $\Gamma $ for which $P$ fails has $p$--mo\-du\-lus zero.

%Recall that $f\in ACP$ if and only if a curve
%$\widetilde{\gamma}=f\circ\gamma$ is rectifiable for a.e. closed
%curves $\gamma$ in $D,$ and a path
%$\widetilde{\gamma}=f(\gamma^0(s))$ is absolutely continuous for
%a.e. closed paths $\gamma$ in $D;$ here $\gamma^0(s)$ denotes the
%normal representation of $\gamma$ defined as above.
A curve $\gamma$ in $D$ is called here a {\it lifting} of a curve
$\widetilde{\gamma}$ in ${\Bbb R}^n$ under $f:D\rightarrow {\Bbb
R}^n$ if $\widetilde{\gamma} = f\circ\gamma.$ Suppose that
$f^{\,-1}(y)$ does not contain a non--degenerate curve for any $y\in
{\Bbb R}^n.$ We say that a mapping $f$ is {\it absolute continuous
on curves in the inverse direction,} abbr. $ACP_p^{-1},$ if for
$p$--a.e. closed curves $\widetilde{\gamma}$ a lifting $\gamma$ of
$\widetilde{\gamma}$ is rectifiable and the corresponding
$f$--re\-pre\-se\-n\-ta\-tion $\gamma^{\,*}$ of $\gamma$ is
absolutely continuous.

\medskip
The result of the paper is the following statement.

\medskip
\begin{theo}\label{th2A} {\sl Let $p\ge 1,$ let a mapping $f:D\rightarrow {\Bbb R}^n$ be
a differentiable a.e., discrete mapping, having $N,$ $N^{-1}$ %$ACP$
and $ACP_p^{-1}$ properties. Suppose that $\Gamma$ is a curve family
in $D,$ that $m$ is a positive integer, and that $f$ winds every
path of $\Gamma$ $m$ times around itself. Then
\begin{equation}\label{eq2A}
M_p(f(\Gamma))\quad\le\quad
\frac{1}{m}\,\cdot\,\int\limits_{D}K_{I,p}(x, f)\cdot
\rho^p(x)\,dm(x)
\end{equation}
for every $\rho\in {\rm adm\,\Gamma}.$ }
\end{theo}

\medskip
Note that the Theorem \ref{th1A} follows from Theorem \ref{th2A} at
$p=n$ as corollary. In fact, every non--constant mapping with
bounded distortion is discrete and has $N$--pro\-per\-ty, see
Theorems 6.2 and 6.3 of Ch. II in \cite{Re} (see also Theorem 4.1
and Proposition 4.14 Ch. I in \cite{Ri}); is differentiable a.e.,
see Lemma 3 in \cite{Va$_1$}; has $N^{\,-1}$--property, see Theorem
8.2 in \cite{BI}; and has $ACP_n^{-1}$--pro\-per\-ty, see Lemma 6 in
\cite{Pol}.

\section{Proof of the main result}\label{sec2}
\setcounter{equation}{0}

A mapping $\varphi:X\rightarrow Y$ between metric spaces $X$ and $Y$
is said to be a {\it Lipschitzian} provided
$$
{\rm dist} \left(\varphi(x_1),\varphi(x_2)\right)\le
M\cdot\text{dist} (x_1,x_2)
$$
for some $M<\infty$ and for all $x_1$ and $x_2\in X.$ The mapping
$\varphi$ is called {\it bi--lipschitz} if, in addition,
$$
M^*\text{dist} \left(x_1,x_2\right)\le\text{dist}
\left(\varphi\left(x_1\right),\varphi\left(x_2\right)\right)
$$ for some $M^*>0$ and for all $x_1$ and $x_2\in X.$ Later on $X$
and $Y$ are subsets of ${\Bbb R}^n$ with the Euclidean distance.

The following proposition can be found in \cite{MRSY$_1$}, see Lemma
3.20, see also Lemma 8.3 Ch. VIII in \cite{MRSY$_2$}.

\medskip
\begin{propo}\label{pr1}{\sl\,
Let $f:D\rightarrow{\Bbb R}^n$ be a differentiable a.e. in $D,$ and
have $N$-- and $N^{-1}$--pro\-per\-ties. Then there is a countable
collection of compact sets $C^*_k\subset D$ such that $m(B_0)=0$
where $B_0=D\setminus\bigcup\limits_{k=1}\limits^{\infty} C^*_k$ and
$f|_{C^*_k}$ is one--to--one and bi--lipschitz for every
$k=1,2,\ldots .$ Moreover, $f$ is differentiable at $C_k^*$ and
$J(x,f)\ne 0.$}
\end{propo}

\medskip
Given a set $E$ in ${\Bbb R}^n$ and a closed curve $\gamma
:\Delta\rightarrow {\Bbb R}^n,$ we identify $\gamma\cap E$ with
$\gamma\left(\Delta\right)\cap E.$
%where $|\gamma |=\gamma (\Delta)$ is the locus of $\gamma .$
If $\gamma$ is rectifiable, then we set
$$
l\left(\gamma\cap E\right) =  m_1(E_ {\gamma}),
$$ where
$E_ {\gamma} = l_{\gamma}\left(\gamma ^{-1}\left(E\right)\right);$
here $l_{\gamma}:\Delta\rightarrow\Delta
_{\gamma}$ as in the previous section. Note that %
$E_ {\gamma} = \gamma^{0\,-1}\left(E\right),$
where $\gamma^0 :\Delta _{\gamma}\rightarrow {\Bbb R}^n$ is the
normal representation of $\gamma $ and
$$l\left(\gamma\cap E\right) = \int\limits_{\Delta} \chi
_E\left(\gamma\left(t\right) \right) |dx|:= \int\limits_{\Delta
_{\gamma}} \chi_{E_\gamma }(s) ds\,.$$

\medskip
The statement mentioned bellow can be found in Chapter IX of
\cite{MRSY$_2$}, see Theorem 9.1.

\medskip
\begin{propo}\label{pr2}{\sl\,
Let $E$ be a set in a domain $D\subset{\Bbb R}^n,$ $n\ge 2.$ Then
$E$ is measurable if and only if $\gamma\cap E$ is measurable for
$p$--a.e. closed curve $\gamma$ in $D.$ Moreover, $m(E)=0$ if and
only if
%\begin{equation}\label{eq2}
$$l(\gamma\cap E)=0$$ %\end{equation}
for $p$--a.e. closed curve $\gamma$ in $D.$}
\end{propo}

\medskip
{\it Proof of the Theorem \ref{th2A}.} Let $B_0$ and $C_{k}^{*}$  be
as in Proposition \ref{pr1}.  Setting by induction $B_1=C_{1}^{*},$
$B_2=C_{2}^{*}\setminus B_1\ldots\,,$
$$B_k\quad=\quad C_{k}^{*}\setminus \bigcup\limits_{l=1}^{k-1}B_l\,,$$
we obtain the countable covering of $D$ consisting of mutually
disjoint Borel sets $B_k,$ $k=1,2,\ldots,$ with $m(B_0)=0.$ By the
construction and $N$--pro\-per\-ty, $m\left(f(B_0)\right)=0.$ Thus,
by Proposition \ref{pr2} $\widetilde{\gamma}^0(s)\not\in f(B_0)$ for
$p$--a.e. curves $\widetilde{\gamma}$ in $f(D),$ where
$\widetilde{\gamma}^{\,0}$ is a normal representation of
$\widetilde{\gamma}.$ Moreover, by $ACP_p^{\,-1}$--pro\-per\-ty, the
$f$--rep\-re\-sen\-ta\-tion $\gamma^{*}$ of a curve $\gamma$ is
rectifiable and absolutely continuous for $p$--a.e. closed curves
$\widetilde{\gamma}$ in $f(D)$ such that
$f\circ\gamma=\widetilde{\gamma}.$

Let $\rho\,\in\,{\rm  \,adm}\,\Gamma$ and
\begin{equation}\label{equa9}
\widetilde{\rho}(y)\quad=\quad\frac{1}{m}\cdot
\chi_{f\left(D\setminus B_0
\right)}(y)\sup\limits_{C}\sum\limits_{x\,\in\,C}\rho^*(x)\,,
\end{equation}
where
$$\rho^*(x)\,=\,\left\{\begin{array}{rr}
\rho(x)/l\left(f^{\prime}(x)\right), &   x\in D\setminus B_0,\\
0,  &  x\in B_0
\end{array}
\right.$$
and $C$ runs over all subsets of $f^{-1}(y)$ in $D\setminus B_0$
such that ${\rm card}\,C\le m.$ Note that
\begin{equation}\label{equa10}
\widetilde{\rho}(y)\quad=\quad\frac{1}{m}\cdot
\sup\sum\limits_{i=1}^s \rho_{k_i}(y)\,,
\end{equation}
where $\sup$ in (\ref{equa10}) is taken over all
$\left\{k_{i_1},\ldots,k_{i_s}\right\}$  such that $k_i\in\,{\Bbb
N},$ $k_i\ne k_j$ if $i\ne j,$ all $s\le m,$  and
$$
\rho_k(y)\,=\left\{\begin{array}{rr}
\,\rho^{*}\left(f_k^{-1}(y)\right), &   y\in f(B_k),\\
0,  &  y\notin f(B_k)
\end{array}
\right. \,.$$
Here $f_k=f|_{B_k},$ $k=1,2,\ldots\,$ is injective and $f(B_k)$ is
Borel. Thus, the function  $\widetilde{\rho}(y)$ is Borel, see e.g.
2.3.2 in \cite{Fe}. %and Theorem I (8.5) in \cite{Sa}.

Given a rectifiable path $\beta$ we denote through $c$ the length of
$\beta,$ $c:=l(\beta),$ and through $\beta^0$ it's normal
representation. Using the definitions of the paths $\Gamma$ and
$f(\Gamma),$ for every curve $\beta\in f(\Gamma)$ we obtain
$$\int\limits_{\beta}\widetilde{\rho}(y)|dy|=
\int\limits_{0}^{c}\widetilde{\rho}(\beta^0(t))dt= $$
$$=\sum\limits_{j=1}^{m}\int\limits_{\frac{c(j-1)}{m}}^{\frac{cj}{m}}\widetilde{\rho}(\beta^0(t))dt=
\sum\limits_{j=1}^{m}\int\limits_{0}^{\frac{c}{m}}\widetilde{\rho}\left(\beta^0\left(t+\frac{c(j-1)}{m}
\right)\right)dt=$$
\begin{equation}\label{eq3A}
=m\int\limits_{0}^{h}\widetilde{\rho}(\beta^0(t))dt,\qquad\,h=c/m\,.
\end{equation}
If $0<t<h,$ then $\alpha^{\,*}(t),$ $\alpha^{\,*}(t+h),\ldots,$
$\alpha^{\,*}(t+(m-1)h)$ are distinct points in
$f^{\,-1}(\beta^0(t)).$ Hence
\begin{equation}\label{eq4A}
\widetilde{\rho}(\beta^0(t))\ge
\frac{1}{m}\sum\limits_{j=0}^{m-1}\rho^{\,*}\left(\alpha^{\,*}(t+jh)\right)
\end{equation}
for $t\in (0, h).$
By Proposition \ref{pr2}, we can consider that $\beta^0(t)\not\in
f(B_0)$ a.e. $t\in [0, c].$ Since $\beta$ is rectifiable, $\beta(t)$
is a differentiable a.e. Besides that, a curve $\alpha^{\,*}$ is
rectifiable and absolutely continuous for $p$--a.e. $\beta,$
moreover, $\alpha^{\,*}(t)\not\in B_0$ for  a.e. $t\in [0, c].$
Thus, the derivatives $f^{\,\prime}\left(\alpha^{\,*}(t)\right)$ and
$\alpha^{\,*\prime}(t)$ exist for a.e. $t\in [0, c].$ Taking into
account the formula of the derivative of the superposition of
functions, and that the modulus of the derivative of the curve by
the natural parameter equals to 1, we have that
$$1=\left|\left(f\circ \alpha^{\,*}\right)^{\,\prime}(t)\right|=
\left|f^{\,\prime}\left(\alpha^{\,*}(t)\right)\alpha^{\,*\prime}(t)\right|=$$
\begin{equation}\label{eq6.6}
=\left|f^{\,\prime}\left(\alpha^{\,*}(t)\right)\cdot\frac{\alpha^{\,*\prime}(t)}{|\alpha^{\,*\prime}(t)|}\right|\cdot
|\alpha^{\,*\prime}(t)|\ge
l\left(f^{\,\prime}\left(\alpha^{\,*}(t)\right)\right)\cdot
|\alpha^{\,*\prime}(t)|\end{equation}
for $p$--a.e. curves $\beta \in f(\Gamma),$ $\beta=f\circ \alpha.$
It follows from (\ref{eq6.6}) that a.e.
\begin{equation}\label{eq6.5}
\rho^{\,*}\left(\alpha^*(t)\right)=\frac{\rho(\alpha^*(t))}
{l\left(f^{\,\prime}\left(\alpha^{\,*}(t)\right)\right)}\ge
\rho(\alpha^*(t))\cdot |\alpha^{\,*\prime}(t)|\,.
\end{equation}

By absolutely continuity of $\alpha^{\,*},$ Theorem 4.1 in
\cite{Va$_2$}, (\ref{eq3A}), (\ref{eq4A}) and (\ref{eq6.5}) we
obtain that
\begin{equation}\label{eq5A}
\int\limits_{\beta}\widetilde{\rho}(y)|dy|\ge
\sum\limits_{j=0}^{m-1}\int\limits_{0}^h \rho(\alpha^*(t+jh))\cdot
|\alpha^{\,*\prime}(t+jh)|dt=\int\limits_{\alpha}\rho(x)|dx|\ge 1\,.
\end{equation}

Thus, $\widetilde{\rho}\,\in\,{\rm }\,\,{\rm
adm}\,f(\Gamma)\setminus \Gamma_0,$ where $M_p(\Gamma_0)=0.$
Consequently
\begin{equation}\label{equa13}
M_p\left(f(\Gamma)\right)\quad
\le\quad\int\limits_{f(D)}\widetilde{\rho}\,^p(y)\,\,dm(y)\,.
\end{equation}
By 2.3.5 for $m=n$ in \cite{Fe}, we obtain that
\begin{equation}\label{equa14}
\int\limits_{B_k}K_{I,
p}(x,\,f)\cdot\rho^p(x)\,\,dm(x)\quad=\quad\int\limits_{f(D)}
\rho^p_k(y)\,dm(y)\,.
\end{equation}
By H\"{o}lder inequality for series,
\begin{equation}\label{equa16}
\left(\frac{1}{m}\cdot\sum\limits_{i=1}^{s}\rho_{k_i}(y)\right)^p\quad\le\quad
\frac{1}{m}\cdot \sum\limits_{i=1}^{s}\,\rho^p_{k_i}(y)
\end{equation}
for each $1 \le s \le m$ and every $\left\{k_1,\ldots,k_s\right\},$
$k_i\in {\Bbb N},$ $k_i\ne k_j,$ if $i\ne j.$

Finally, by the Lebesgue positive convergence theorem, see Theorem
12.3 in $\S\,12.$ I in \cite{Sa}, we conclude from
(\ref{equa13})--(\ref{equa16}) that
%
%\begin{equation}\label{equa15}
$$\frac{1}{m}\cdot\int\limits_{D}K_{I,
p}(x,\,f)\cdot\rho^p(x)\,\,dm(x)\quad
=\quad\frac{1}{m}\cdot\int\limits_{f\left(D\right)}\,\sum\limits_{k=1}^{\infty}
\rho_k^p(y)\,dm(y)\quad\ge$$
%\end{equation}
%
$$\ge\quad\frac{1}{m}\cdot\int\limits_{f\left(D\right)}
\sup\limits_{\left\{k_1,\ldots,k_s\right\},\, k_i\in {\Bbb N},\atop
k_i\ne k_j \,{\rm if}\, i\ne j}\sum\limits_{i=1}^s
\rho^p_{k_i}(y)\,dm(y)\quad\ge\quad
\int\limits_{f\left(D\right)}\,\widetilde{\rho}^{\,p}(y)\,dm(y)\quad\ge\quad
M_p(f(\Gamma))\,.$$
The proof is complete. $\Box$

\begin{rem}\label{rem1}
The above investigations are closely related with the so--called
mappings with finite length distortion, introduced by O. Martio
together with V. Ryazanov, U. Srebro and E. Yakubov, see
\cite{MRSY$_1$}--\cite{MRSY$_2$}, see also works \cite{BGMV} and
\cite{KO}.
\end{rem}

%=================Список литературы====================
%\end{fulltext}

{\bf \noindent Evgeny A. Sevost'yanov} \\
Institute of Applied Mathematics and Mechanics,\\
National Academy of Sciences of Ukraine, \\
74 Roze Luxemburg str., 83114 Donetsk, UKRAINE \\
Phone: +38 -- (062) -- 3110145, \\
Email: esevostyanov2009@mail.ru
\end{document}